\documentclass[12pt]{iopart}

\usepackage{graphicx}
\usepackage{multicol}
\usepackage{dsfont}
\usepackage{mathrsfs}
\usepackage{amssymb}
\usepackage{upgreek}
\usepackage{esint}
\newtheorem{remark}{Remark}
\newtheorem{theorem}{Theorem}

\begin{document}

\title[]{On the Computation of Hierarchical Control results for One-Dimensional Transmission Line}

\author{P. P. de Carvalho \ \&  \ O. P. de S\'a Neto}

\address{Universidade Estadual do Piau\'i - UESPI}
\ead{pitagorascarvalho@gmail.com  \\ \ \ \ \ \ \ \ olimpioqedc@gmail.com}

\vspace{10pt}
\begin{indented}
\item[]  2021
\end{indented}

\begin{abstract}
In this paper, motivated by a physics problem, we investigate some numerical and computational aspects for the problem of hierarchical controllability in a one-dimensional wave equation in domains with a moving boundary. Some controls act in part of the boundary and define a strategy of equilibrium between them, considering a leader control and a follower. Thus, we introduced the concept of hierarchical control to solve the problem and mapped the Stackelberg Strategy between these controls.
A total discretization of the problem is presented for a numerical evaluation in spaces of finite dimension, an algorithm for evaluation of the problem is presented as the combination of finite element method (FEM) and finite difference method (FDM). The algorithm efficiency and computational results are illustrated for some experiments using the software \textsc{FreeFem\raisebox{.4ex}{\tiny\bf ++}}.
\end{abstract}

\vspace{2pc}
\noindent{\it Keywords}: Transmition Line Ressonator, Wave Equation, Moving boundary, Hierarchical control, Stackelberg strategy, Finite element method.
\maketitle

\section{Introduction}

On several occasions, controlling a problem may involve more than one agent (control). For such situations, we can define a strategy that indicates the desired behavior.
This paper deals with the numerical solution of a controllability problem for the wave equation
through a hierarchy of controls in boundary. More precisely, we have chosen the so called \textit{Stackelberg-Nash} method,
that can be briefly described as follows:
\begin{itemize}
  \item We have control of two kinds: leaders and followers.
  \item  We associate to leader a \textit{Nash} equilibrium, that corresponds to a noncooperative
multiple-objective optimal control problem.
  \item Then, we choose the leader among the set of controls by minimizing a suitable functional.
\end{itemize}

Initially, in game theory a player is a strategic decision-maker within the context of the game.
And the game is characterized by any set of circumstances that have an outcome depends on the actions of two or more decision-makers (players).
In a hierarchical game, that is, in which all players make their decisions based on a decision by a leading player, and a result is achieved for all players, this is called the equilibria position. In our case, there will be no cooperation in decision making between players. In other words, fixed a leader we will dedicate ourselves to the study of equilibria where there is a leader and the other players adopt Nash equilibrium in the equations.
The process in the problems above is a combination of strategies and is called Stackelberg-Nash strategy.
For more details in noncooperative optimization strategy proposed by Nash (see \cite{N}) and the Stackelberg hierarchical-cooperative strategy (see \cite{stalck}).

Some numerical and computational results involving \textit{Nash} equilibrium we can found in \cite{Cara-Pita1}, \cite{RA1} and \cite{RA2}, the \textit{Stackelberg-Nash} equilibrium  in \cite{Cara-Pita3}. For the algorithm construction, we adapt the ideas contained in \cite{Cara-Pita2} in that the authors work in numerical viewpoint the \textit{Nash equilibrium} for the wave equations, but the controls domains acting in subregions of the domain.

The structure of the article is given as follows: In Section \ref{PS} and \ref{sec2}, we present a physical motivation for the problem and the system of control respectively. Sections \ref{sec3} and \ref{sec4}, are devoted to some comments for the existence and uniqueness of \textit{Nash equilibrium} and present the approximate controllability with respect to the leader control.  Section \ref{sec5} we present the optimality system for the leader control, the principal results for the strategy of \textit{Nash} for the linear system obtained in \cite{Je1000}. Section \ref{Sec6} we leave it reserved to full discretization and presentation of the algorithm used to solve the problem. Section \ref{Sec7} concentrates tables and numerical experiments resulting from the data simulation presented in Section \ref{Sec6}.
Finally, in Section \ref{Sec8} some comments and possible advances are added.

\section{Physical Systems}
\label{PS}

\begin{figure}[h]
\begin{center}
\includegraphics[scale=0.3]{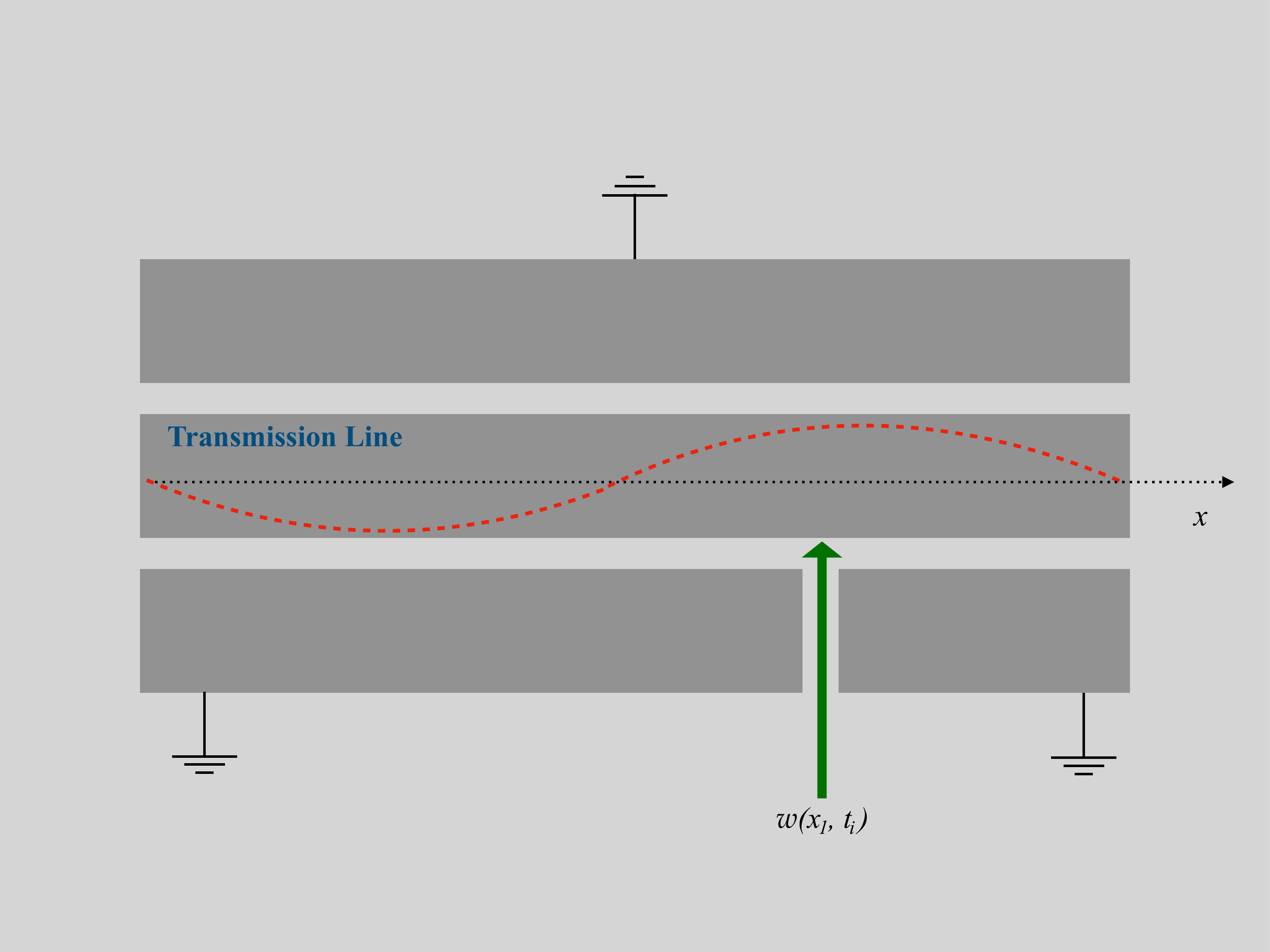}
\caption{Modes of load density vibrations in a transmission line in schematic model for spatially located control in time did not continue. }
\label{fig}
\end{center}
\end{figure}

When we transmit a microwave signal through a $l$ length transmission line, if the wavelength is much greater than the cross-sectional dimension of the line, the loads on the transmission line can be considered as if they were moving in a single dimension, figure (\ref{fig}). The \textit{n} radiation modes behaved in this transmission line can be modeled by a set of discrete and infinitesimal LC elements known as concentrated circuit elements \textit{(lumped circuit)}\cite{DP}.

The Lagrangean in the system is:

\begin{eqnarray}
\displaystyle \mathcal{L} &=& \sum_{n}\bigg[\frac{li_{n}^{2}}{2}-\frac{q_{n}^{2}}{2c}\bigg],
\label{2}
\end{eqnarray}
where $c$ is the capacitance and $l$ is the auto inductance of $n$ - this is the mode of the transmission line. In this case, the temporal variation of the load at the $n$ node of the circuit is given by $\dot{q_{n}}=i_{n-1}-i_{n}$ and $i_{n}=-\sum^{n}_{m=1}\dot{q_{m}}$ is the current at the node. Substituting in Lagrangean (\ref{2}), we have
\begin{eqnarray}
\displaystyle \mathcal{L} &=& \sum_{n}\Bigg[ \frac{l(\sum_{m=1}^{n}\dot{q}_{m})^{2}}{2}-\frac{q_{n}^{2}}{2c}\Bigg].
\label{3}
\end{eqnarray}
As usual in this type of system we will make use of the infinitesimal nature of these elements (the degrees of freedom of the system) to take the equation (\ref{3}) into the continuum. We define the variable
\begin{equation}
\displaystyle u(x,t)=\int_{-\frac{L}{2}}^{x} dx^{'} q(x^{'},t)
\label{4}
\end{equation}
where $q(x)$ is the linear density of charge. Making the substitutions:
$$\displaystyle \sum^{m=1}_{n}q_{m}(t)\rightarrow u(x,t), \ \ \ \ \ q_{n}(t)\rightarrow q(x,t)=\frac{\partial u}{\partial x},$$
one-dimensional Lagrangean density is written\begin{eqnarray}
\mathcal{L} &=& \frac{l\dot{u}^2}{2}-\frac{1}{2c}\left[\frac{\partial\Gamma}{\partial x}\right]^2.
\label{5}
\end{eqnarray}
Here $c$ and $l$ are transformed into linear capacitance density and transmission line inductance, respectively. Applying Euler-Lagrange to (\ref{5}), we obtain
\begin{eqnarray}
\frac{1}{c} \frac{\partial^2 u}{\partial x^2}-l \frac{\partial^2 u}{\partial t^2} &=& 0,
\label{6}
\end{eqnarray}
where $1/\sqrt{lc}$ is the velocity of the wave propagation. For our present problem, we will consider the equation of the wave with dimensionless velocity $1/\sqrt{lc}=1$, and simplifications of annotations for
\begin{eqnarray*}
\frac{\partial^2 u}{\partial x^2} &\equiv u_{xx},
\end{eqnarray*}
and
\begin{eqnarray*}
\frac{\partial^2 u}{\partial t^2} &\equiv u_{tt}.
\end{eqnarray*}

\section{Statement of the problem}\label{sec2}
\noindent

Initially, we consider the non-cylindrical domain as constructed in \cite{Je1000}:
$$\displaystyle \widehat{Q}= \left\{ (x,t) \in {\mathds{R}}^2; \; 0 < x <  \alpha _{k(t)},\; \; t \in (0,T) , \ \ \ T> 0\right \},$$
with
$$ \alpha_{k}(t) = 1 + kt, \;\;\;\;\;\; \; 0 < k < 1 \ ,$$
 the lateral boundary defined by $\displaystyle \widehat{\Sigma} =\widehat{\Sigma}_0
\cup \widehat{\Sigma}_0^*$, where
$$\displaystyle \widehat{\Sigma}_{0} = \{(0, t);\; t \in (0,T) \} \;\;\;\;\mbox{ and } \;\;\;\; \widehat{\Sigma}_{0}^* =\widehat{\Sigma} \backslash \widehat{\Sigma}_0= \{(\alpha_k(t),t);\;  t \in (0,T) \}.$$
Consider $\Omega_t$ and $\Omega_0$ the intervals $\displaystyle (0, \alpha_k(t))$
and $\displaystyle (0, 1)$ respectively, and the following system in the domain
$\widehat{Q}$:
\begin{equation} \label{eq1.3}
\left\{
\begin{array}{l}
\displaystyle u_{tt} - u_{xx} = 0 \ \ \mbox{ in } \ \ \widehat{Q},\\[11pt]
\displaystyle u(x,t)\Big{|}_{\widehat{\Sigma}_{0}} = \widetilde{w}(t) \ \ \mbox{and} \ \ u(x,t)\Big{|}_{\widehat{\Sigma}_0^*} = 0,\\[12pt]
\displaystyle u(x,0) = u_0(x), \;\; u_{t}(x,0) = u_1(x) \ \mbox{ in }\; \Omega_0,
\end{array}
\right.
\end{equation}
with $u$ the state, $\widetilde{w}$ the control and $( u_0(x),
u_1(x)) \in L^2(0,1) \times H^{-1}(0,1)$.

The problem (\ref{eq1.3}) models the motion of a string where an endpoint is fixed and the other one is moving and the constant $k$ is called the speed of the moving endpoint.

Consider
\begin{equation}\label{decomp0}
\displaystyle \widehat{\Sigma}_0=\displaystyle \widehat{\Sigma}_1 \cup \displaystyle \widehat{\Sigma}_2, \ \ \ \mbox{with} \ \ \ \displaystyle \widehat{\Sigma}_1 \cap \displaystyle \widehat{\Sigma}_2 = \emptyset
\end{equation}
 and
\begin{equation} \label{decomp}
\displaystyle \widetilde{w}=\{\widetilde{w}_1, \widetilde{w}_2\}, \;\; \widetilde{w}_i=\mbox{control function in } \; L^2(\widehat{\Sigma}_i), \;i=1,2.
\end{equation}

We can also write
\begin{equation} \label{decomp 2.A}
\displaystyle \widetilde{w}= \widetilde{w} _1 + \widetilde {w}_2, \; \mbox{ with } \; \displaystyle {\widehat \Sigma}_0=\displaystyle {\widehat \Sigma}_1 = \displaystyle {\widehat \Sigma}_2.
\end{equation}

Can be rewritten the system (\ref{eq1.3}) as follows:
\begin{equation} \label{eq1.3.1}
\left\{
\begin{array}{l}
\displaystyle u_{tt}- u_{xx} = 0 \ \ \mbox{ in } \ \ \widehat{Q},\\[5pt]
\displaystyle u(x,t)\Big{|}_{\widehat{\Sigma}_{1}} = \widetilde{w}_1(t), \ \ \ u(x,t)\Big{|}_{\widehat{\Sigma}_{2}} = \widetilde{w}_2(t)
 \ \ \mbox{ and } \ \ \displaystyle u(x,t)\Big{|}_{\widehat{\Sigma}\backslash \widehat{\Sigma}_0} = 0,
\\[10pt]
\displaystyle u(x,0) = u_0(x), \;\; u_t(x,0) = u_1(x) \ \mbox{ in }\; \Omega_0.
\end{array}
\right.
\end{equation}
Consider $\widetilde{w}_1$ as being the ``main'' control (the leader), $\widetilde{w}_2$ as the follower, in \textit{Stackelberg} terminology and $u=u(x,t)$ the solution of (\ref{eq1.3.1}).
We will also introduce the (secondary) functional
\begin{equation}\label{sfn}
\displaystyle \widetilde{J}_2(\widetilde{w}_1, \widetilde{w}_2) =
 \frac{1}{2} \iint_{\widehat{Q}} \left(u(\widetilde{w}_1, \widetilde{w}_2)- {u}_2 \right)^2 \,dx\,dt + \frac{\sigma}{2} \int_{\widehat{\Sigma}_2} \widetilde{w}_2^2\;d\widehat{\Sigma},
\end{equation}
and the (main) functional
\begin{equation}\label{mfn}
\displaystyle \widetilde{J}(\widetilde{w}_1) =\frac{1}{2}\int_{\widehat{\Sigma}_1} \widetilde{w}_1^2\;d\widehat{\Sigma},
\end{equation}
where ${\sigma}>0$ is a constant and ${u}_2$ is a given
function in $L^2(\widehat{Q}).$

\begin{remark}\label{rsol} As  in \textup{\cite{Mi}}, we can prove that for each $u_0 \in L^2(0,1)$, $u_1 \in H^{-1}(0,1)$ and $\displaystyle \widetilde {w}_i \in
L^2(\widehat{\Sigma}_i), \;i=1,2,$ there exists exactly one  solution $\displaystyle u$ to
(\ref{eq1.3.1}) in the sense of a transposition, in particular, the cost
functionals $\displaystyle \widetilde{J}_2$ and $\displaystyle \widetilde {J}$ are well defined.
\end{remark}

The \textit{Stackelberg-Nash} strategy:  Thus, if the leader $ \widetilde{w}_1$ makes a choice, then the follower $ \widetilde{w}_2$ makes also a choice, depending on  $\widetilde{w}_1$, which minimizes the cost
 $\widetilde{J}_2$, that is,
\begin{equation}\label{son}
\displaystyle \widetilde{J}_2(\widetilde{w}_1, \widetilde{w}_2)= \inf_{\widehat{w}_2 \in L^2(\widehat{\Sigma}_2)} \widetilde{J}_2(\widetilde{w}_1, \widehat{w}_2).
\end{equation}

%%%%%%%%%%%%%%%%%%%%%%%%%%%%%%%%%%%%%%%%%%%%%%%%%%%%%%%%%%%%%%%%%%%%%%%%%%%%%%%%%%%%%%%%%%%%%
\section{Nash equilibrium}\label{sec3}

\noindent

In this section, fixed any leader control $\displaystyle \ w_1 \in L^2(\widehat{\Sigma}_1)$ we
determine the existence and uniqueness of solutions to the problem
\begin{equation} \label{eq3.10}
\begin{array}{l}
\displaystyle\inf_{\widetilde{w}_2 \in L^2(\widehat{\Sigma}_2)}J_2(\widetilde{w}_1,\widetilde{w}_2),
\end{array}
\end{equation}
and a characterization of this solution in terms of an adjoint system.

In fact, this is a classical type problem in the control of distributed systems (cf.
J. -L. Lions \cite{L3}). It admits an unique solution
\begin{equation} \label{eq3.11}
\displaystyle \widetilde{w}_2 = \mathfrak{F}(\widetilde{w}_1).
\end{equation}
The Euler - Lagrange equation for problem  (\ref{eq3.10}) is given by
\begin{equation} \label{eq3.21}
\int_{0}^{T}\int_{\Omega_t}(u - {u}_2)\widehat{u}dxdt + {\sigma}\int_{\widehat \Sigma_2}\widetilde{w}_2\widehat{w}_2d\ {\widehat \Sigma} = 0, \;\;\forall\, \widehat{w}_2 \in L^2(\widehat \Sigma_2),
\end{equation}
where $\widehat{u}$ is solution of the following system
\begin{equation} \label{eq3.22}
\left\{
\begin{array}{l}
\displaystyle \widehat u_{tt} - \widehat u_{xx} = 0 \ \ \mbox{ in } \ \ \widehat{Q},\\ [7pt]
\displaystyle \widehat{u}\Big{|}_{\widehat{\Sigma}_1} = 0 , \ \ \widehat{u}\Big{|}_{\widehat{\Sigma}_2}=\widehat{w}_2 \ \ \mbox{and} \ \ \widehat{u}\Big{|}_{\left({\widehat\Sigma_1} \cup {\widehat \Sigma_2}\right)} = 0 \\[11pt]
\displaystyle \widehat{u}(x,0) = 0, \; \widehat{u}_{t}(x,0) = 0, \;\; x \in \Omega_t.
\end{array}
\right.
\end{equation}
In order to express (\ref{eq3.21}) in a convenient form, we introduce the adjoint to (\ref{eq3.22})
state defined by
\begin{equation}\label{sac}
\left\{
\begin{array}{l}
\displaystyle p_{tt} - p_{xx} = u - {u}_2  \ \mbox{ in } \ \ \widehat{Q}, \\[5pt]\displaystyle
p(T) = p_{t}(T) = 0, \;\; x \in \Omega_t, \\[10pt]\displaystyle
p = 0 \ \mbox{ on } \ \widehat \Sigma.
\end{array}
\right.
\end{equation}
Multiplying (\ref{sac}) by $\widehat{u}$ and integrating by parts, we find
\begin{equation} \label{eq3.33}
\int_{0}^{T}\int_{\Omega_t} (u - {u}_2)\widehat{u}\,dx\,dt + \int_{\widehat \Sigma_2} p_x\,\widehat{w}_2\,d{\widehat \Sigma} = 0,
\end{equation}
so that (\ref{eq3.21}) becomes
\begin{equation}\label{ci}
\displaystyle p_x= {\sigma} \,\widetilde{w}_2 \ \ \mbox{ on } \ \ \widehat \Sigma_2.
\end{equation}
We summarize these results in the following theorem.
\begin{theorem}\label{teN} For each $\displaystyle \widetilde {w}_1 \in L^2(\Sigma_1)$ there exists a unique \textit{Nash equilibrium} $\displaystyle \widetilde{w}_2$ in the sense of (\ref{son}). Moreover, the follower $\displaystyle \widetilde{w}_2$ is given by
\begin{equation}\label{cseg}
\displaystyle \displaystyle \widetilde{w}_2 = \mathfrak{F}(\widetilde {w}_1)=\frac{1}{{\sigma}}\,\;p_x\;\;\mbox{ on }\;\; \widehat{\Sigma}_2,
\end{equation}
where $\displaystyle \{ v,p \}$ is the unique solution of (the optimality system)
\begin{equation} \label{eq3.37}
\left\{
\begin{array}{l}
\displaystyle u_{tt} - u_{xx} = 0 \ \mbox{ in } \;\; \widehat{Q}, \\
\displaystyle p_{tt} - p_{xx} = u - {u}_2 \ \mbox{ in } \;\; \widehat{Q},\\[5pt]
\displaystyle u \Big{|}_{\widehat{\Sigma}_1} = \widetilde {w}_1, \ \ \displaystyle u \Big{|}_{\widehat{\Sigma}_2} =  \frac{1}{\widetilde{\sigma}}\;\,p_x
\ \ \mbox{and} \ \ \displaystyle u \Big{|}_{\widehat{\Sigma} \backslash \widehat{\Sigma}_0} = 0,
\\[7pt]
p = 0 \ \mbox{ on } \ \widehat{\Sigma}, \\[5pt]
u(0) = u_{t}(0) = 0, \\[5pt]
p(T) = p_{t}(T) = 0, \;\; x \in \Omega_t.
\end{array}
\right.
\end{equation}
Of course, $\displaystyle \{ u,p \}$ depends on $\widetilde{w}_1$:
\begin{equation}\label{cdep}
\displaystyle \{ u,p \} = \{ u(\widetilde{w}_1),p(\widetilde {w}_1) \}.
\end{equation}
\end{theorem}

\section{On the approximate controllability}\label{sec4}

\noindent

Since we have proved the existence, uniqueness and characterization of the
follower $\displaystyle \widetilde{w}_2$, the leader $\displaystyle \widetilde{w}_1$  now wants  that the solutions $u$
and $u'$, evaluated at time $t=T$, to  be as close as possible to $\displaystyle (u^0,
u^1)$. This will be possible if the system (\ref{eq3.37})  is approximately
controllable.

 We are looking for
\begin{equation} \label{inf1}
\begin{array}{l}
\displaystyle\inf\, \frac{1}{2\,}\,\int_{\widehat \Sigma_1} \widetilde{w}_{1}^{2}\,d{\widehat \Sigma},
\end{array}
\end{equation}
where $\displaystyle \widetilde{w}_1$ is subject to
\begin{equation} \label{subj1}
\begin{array}{l}
\displaystyle \left(u(T;{\widetilde{w}_1}), u'(T; {\widetilde{w}_1})\right)\in B_{L^2(\Omega_t)}(u^0,\rho_0) \times B_{H^{-1}(\Omega_t)}(u^1,\rho_1),
\end{array}
\end{equation}
assuming that  $w_1$ exists, $\rho_0$ and  $\rho_1$ being positive numbers arbitrarily
small  and $\{u^0, u^1\} \in L^2(\Omega_t) \times H^{-1}(\Omega_t)$.

As in \cite{Je1000}, we assume that
\begin{equation}\label{hT}
\displaystyle T >  \frac{e^{\frac{2k(1 + k)}{(1 - k)^3}} - 1}{k}
\end{equation}
and
\begin{equation}\label{hT10}
0 < k <  1.
\end{equation}

\begin{theorem}{Assume that (\ref{hT}) and (\ref{hT10}) hold. Let us consider $\displaystyle \widetilde{w}_1 \in L^2(\widehat \Sigma_1)$ and $\displaystyle \widetilde{w}_2$ a \textit{Nash equilibrium} in the sense (\ref{son}).
Then $$\displaystyle \left(u(T), u'(T)\right)=\left(u(., T, {\widetilde{w}_1}, \widetilde{w}_2), v'(., T, {\widetilde{w}_1},
\widetilde{w}_2)\right) \, ,$$ where $\displaystyle u$ solves the system (\ref{eq3.37}), generates a dense
subset of $\displaystyle L^2(\Omega_t)\times H^{-1}(\Omega_t)$.}
\end{theorem}\label{AC}
\begin{remark}
As can be seen in \cite{Je1000}, the income statement above is done using the decomposition of the solutions in (\ref{eq3.37})
\begin{equation} \label{eq3.39}
\left\{
\begin{array}{l}
u = u_0 + g,\\
p = p_0 + q,
\end{array}
\right.
\end{equation}
where $u_0$, $p_0$,$g$ and $q$ are particular solutions for this system.
New systems for $g$ and $q$ are obtained, and the author
consider the following ``adjoint systems" for $g$ and $q$ respectively:
\begin{equation} \label{eq3.45}
\left\{
\begin{array}{l}
\varphi_{tt} - \varphi_{xx} =\displaystyle\psi \ \mbox{ in } \ \widehat Q, \\[5pt]\displaystyle
\varphi = 0 \ \mbox{ on } \ \widehat \Sigma,\\[5pt]\displaystyle
\varphi(T) = 0, \ \varphi_{t}(T) = 0, \;\; x \in \Omega_t,
\end{array}
\right.
\end{equation}
and
\begin{equation} \label{eq3.46}
\left\{
\begin{array}{l}
\psi_{tt} - \psi_{xx} = 0 \ \mbox{ in } \ \widehat Q,\\[5pt]\displaystyle
\psi \Big{|}_{\widehat \Sigma_1} = 0, \ \ \ \psi \Big{|}_{\widehat \Sigma_2} = \displaystyle \frac{1}{\sigma}\,\varphi_{x} \ \ \ \mbox{and} \ \ \ \psi \Big{|}_{ \widehat\Sigma \backslash \widehat \Sigma_0} = 0,
\\[10pt]\displaystyle
\psi(0) = \psi_{t}(0) = 0, \;\; x \in \Omega_t.
\end{array}
\right.
\end{equation}
\end{remark}

%%%%%%%%%%%%%%%%%%%%%%%%%%%%%%%%%%%%%%%%%%%%%%%%%%%%%%%%%%%%%%%%%%%%%%%%%%55
%%%%%%%%%%%%%%%%%%%%%%%%%%%%%%%%%%%%%%%%%%%%%%%%%%%%%%%%%%%%%%%%%%%%%%%%%%

\section{Optimality systems and main results}\label{sec5}
\noindent

Thanks to the results obtained in preceding sections, we can achieve for  each
$\displaystyle \widetilde{w}_1$, the \textit{Nash equilibrium} $\displaystyle \widetilde{w}_2$ associated to solution $\displaystyle u$
of (\ref{eq1.3.1}).

Let us consider
\begin{equation} \label{eq3.7cil.1}
\displaystyle \inf_{\widetilde{w}_1\in \mathcal{U}_{ad}} J(\widetilde{w}_1),
\end{equation}
where $\displaystyle \mathcal{U}_{ad}$ is the set of admissible controls
\begin{equation}\label{admcon}
\displaystyle \mathcal{U}_{ad}=\{\widetilde{w}_1\in L^2({\widehat \Sigma_1}); \; u \mbox{ solution of } (\ref{eq1.3.1}) \mbox{ satisfying } (\ref{subj1})\}.
\end{equation}

Again as in \cite {Je1000}, the following result holds:
\begin{theorem} \label{teor3.6} Assume the hypotheses (\ref{decomp 2.A}), (\ref{hT}) and (\ref{hT10})  are satisfied. Then for $\{f^0,f^1\}$
in $\displaystyle H_0^1(\Omega_t) \times L^2(\Omega_t)$ we uniquely define $\{\varphi, \psi, u, p \}$
by
\begin{equation} \label{eq3.139}
\left\{
\begin{array}{l}
\displaystyle \varphi_{tt} - \varphi_{xx} = \psi \ \ \mbox{in} \ \ \widehat Q, \\[3pt]\displaystyle
\displaystyle \psi_{tt} - \psi_{xx}  = 0 \ \ \mbox{in} \ \ \widehat Q, \\[3pt]\displaystyle
\displaystyle u_{tt} - u_{xx} = 0 \ \ \mbox{in} \ \ \widehat Q, \\[3pt]\displaystyle
\displaystyle p_{tt} - p_{xx}  = u - {u}_2 \ \ \mbox{in} \ \ \widehat Q, \\[3pt]\displaystyle
\displaystyle \varphi = 0 \ \ \ \mbox{on} \ \ \ \widehat \Sigma, \\[3pt]\displaystyle
\displaystyle \psi \Big{|}_{\widehat \Sigma_1} = 0, \ \ \ \displaystyle \psi \Big{|}_{\widehat \Sigma_2} = \displaystyle \frac{1}{\sigma}\;\varphi_{x} \ \ \
 \mbox{and} \ \ \ \psi \Big{|}_{\widehat \Sigma \backslash \widehat \Sigma_0} =  0, \\[1pt]\displaystyle
\\
\displaystyle u \Big{|}_{\widehat \Sigma_1} =  - \varphi_{x}, \ \ \ \displaystyle u \Big{|}_{\widehat \Sigma_2} = \displaystyle \frac{1}{ \sigma}\;p_{x} \ \ \ \mbox{and} \ \ \ \displaystyle u \Big{|}_{\widehat \Sigma \backslash \widehat \Sigma_0}=0 ,\\[3pt]\displaystyle
\displaystyle p = 0 \ \ \ \mbox{on} \ \ \ \widehat \Sigma, \\[3pt]\displaystyle
\displaystyle \varphi(.,T) = 0,\, \varphi_{t}(.,T) = 0 \ \ \mbox{in} \ \  \Omega_t, \\[3pt]\displaystyle
\displaystyle u(0) = u_{t}(0) = 0 \ \ \mbox{in} \ \   \Omega_t, \\[3pt]\displaystyle
\displaystyle \psi(0) = \psi_{t}(0) = 0 \ \ \mbox{in} \ \   \Omega_t, \\[3pt]\displaystyle
\displaystyle p(T) = p_{t}(T) = 0 \ \ \mbox{in} \ \  \Omega_t.
\end{array}
\right.
\end{equation}
As in \cite{Je1000}, the optimal leader is given by
\begin{equation*}
\widetilde{w}_1 = - \varphi_{x} \ \ \mbox{on} \ \  \widehat \Sigma_1,
\end{equation*}
where $\varphi$ corresponds to the solution of first equation in the system (\ref{eq3.139}).
\end{theorem}

\section{Full discretization and Algorithm}\label{Sec6}

%\subsection{A Fixed Point Algorithm for Solving \eqref{eq3.139}} \label{Sec4.1}
\noindent

In the sequel we employ a methodology combining finite
differences for the time discretization, finite elements for the space approximation, and a fixed point
algorithm for the iterative solution of the discrete control problem for (\ref{eq3.139}), using ideas similar to those developed in \cite{RA1}, \cite{Cara-Pita1} and \cite{Cara-Pita2}.
%\vspace{.2cm}

Initially, introduce the notation $$V_i := L^2(\widehat \Sigma_i) \ \ \mbox{and} \ \ V = V_1 \times V_2 \, , $$ where $V:= L^2(\widehat \Sigma_0)$.

As consequence of anterior results we have that: for all ${\sigma} > 0$ (sufficiently large), exists an unique equilibrium $(\tilde{w}_1, \, \tilde{w}_2) \in V$ for the functionals $\tilde{J}_1$ and $\tilde{J}_2$ , satisfying the Theorem \ref{teor3.6}.
%and $\mathcal{H} := \mathcal{H}_{1} \times \mathcal{H}_{2}$ and let us consider the operator $L_{i} \in \mathcal{L}(\mathcal{H}_{i} \ ; L^2(Q))$ with $L_{i}v_{i} = z_{i}$, where $z_{i}$ is the solution to the system

\subsection{Reduction to Finite Dimension}\label{Sec3.3}
\noindent

In what follows, we will describe some approximate spaces and schemes in the next section.

The reduction of (\ref{eq3.139}) to finite dimension must be performed in two steps:
\noindent

  $\bullet$ {\sc Step 1: Approximation in time.}  We consider the time discretization step $ \Delta t$, defined by $ \Delta t = T/M $, where $M$ is a large positive integer. Then, if we set $ t^m := m \Delta t$, we have
$$ 0 < t^1 < t^2 < \cdot \cdot \cdot < t^{M} = T .$$
 Now, we approximate $V_{1}$ and $V_{2}$ respectively by
 $${V}_{1}^{\Delta t} := L^2(\widehat \Sigma_1)^M  \ \ \ \ \ \ \mbox{and} \ \ \ \ \ \ {V}_{2}^{\Delta t} := L^2(\widehat \Sigma_2)^M .$$

Accordingly, we can interpret the elements of \, $V^{\Delta t}_{i}$ as controls in $V_{i}$ that are piecewise constant in time.

\noindent

$\bullet$ {\sc Step 2: Approximation in space.} From now on, we will establish to fix ideas that $\Omega_t$ is a subdomain of $\mathds{R}$. We will also assume that the $\widehat \Sigma$ is the total boundary  and $\widehat \Sigma_i$ are boundary segments which denote the domains of controls $\tilde{w}_i$ (with $i=1,2$). We introduce a triangulation $\mathcal{T}_{h}$ of $\Omega_t$, where we assume that $h$ is the longest length of the edges of the triangles of  $\mathcal{T}_{h}$ .
Next, we approximate the set of solutions
$$W= \{w \in L^{\infty}(0,T; H^{1}_{0}(\Omega_t)) \ : \ w_{t} \in L^{\infty}(0,T;L^2(\Omega_t))\}$$
by  $W^{\Delta t}_{h},$ where
 $$ W^{\Delta t}_{h} := (W_{h})^M , \, \ \ \ \  W_{h} := \bigl\{ z \in \mathcal{C}^{0}({\Omega_t}) \, : \ z \big{|}_{K} \in \mathds{P}_{1}(K) \, \ \forall  K \in\mathcal{T}_{h} \bigr \}  $$
and $\mathds{P}_{1}(K)$ is the space of the polynomial functions of degree $\leq 1$; thus, $\dim(\mathds{P}_{1}(K)) = 3$ \ and \ $\dim(W_{h}) = N_{h},$  \ where $N_{h}$ is the number of vertices of $\mathcal{T}_{h}$.
In this second step, we first approximate $V^{\Delta t}_{i}$ by $V^{\Delta t}_{i,h}$, defined as follows:
 $$ V^{\Delta t}_{i,h} = (V_{i,h})^{M} , \ \; \; \ {V}_{i,h} := \bigl\{ z \in \mathcal{C}^{0}({\widehat \Sigma}_{i}) \, : \ z \big{|}_{K} \in \mathbb{P}_{1}(K) \; \ \forall  K \in \widehat \Sigma_{i} \bigr \};$$
then, we set ${V}^{\Delta t}_{h} := {V}^{\Delta t}_{1,h} \times {V}^{\Delta t}_{2,h}$.
Finally, we consider the finite-dimensional version of $W$, determined by
 $$W^{\Delta t}_{h,0} = (W_{h,0})^N , \ \ \ \  W_{h,0} = \bigl \{ z \in W_{h} : z \big{|}_{\Gamma_{1}} = 0 \bigr\}\, .$$
In (\ref{eq3.139}) the state equation (in $u$ and $\psi$) and the adjoint systems (in $p$ and $\varphi$) can be approximated in time and space incorporating
(for instance) implicit Euler finite differences in time and spatial $P_1$-Lagrange finite element techniques. That allows us to compute a state $u^{\Delta t}_{h}$ and two adjoint states $p^{\Delta t}_{h}$ and $\varphi^{\Delta t}_{h}$ for each control pair $\widetilde{w} = (\widetilde{w}_1, \widetilde{w}_2) \in {V}^{\Delta t}_{h}$.

In accordance with  those definitions, we can approximate the problem to obtain a pair control $(\widetilde{w}_1, \, \widetilde{w}_2) \in V$ by a finite dimensional problem:
\begin{equation}\label{3.10a}
\left\{
\begin{array}
[c]{ll}%
\vspace{.3cm}
 \displaystyle\frac{\partial \widetilde{J}^{\Delta t}_{1,h}}{\partial \widetilde{w}_{1}}(\widetilde{w}_{1}, \widetilde{w}_{2}) = 0 , \\
\vspace{.08cm}
 \displaystyle\frac{\partial \widetilde{J}^{\Delta t}_{2,h}}{\partial \widetilde{w}_{2}}(\widetilde{w}_{1}, \widetilde{w}_{2}) = 0 ,
\end{array}
\right.
   \end{equation}
where the $\widetilde{J}^{\Delta t}_{i,h}$ are the finite-dimensional versions of the $\widetilde{J}_i$ induced by time and space approximations.

\subsection{Fixed--Point Method for the Discretized Linear Problem}\label{Sec3.4}

Now, we can solve the approximate formulation (\ref{3.10a}) for a the equivalent problem using the fixed--point algorithm as follows:

\vspace{.15cm}
\begin{description}
  \item[ALGORITHM:]
  \item[ a)] Choose ${\widetilde{w}}_0 := ({\widetilde{w}}_{1,0} , \, {\widetilde{w}}_{2,0} ) \in V^{\Delta t}_{h}$ (where ${\widetilde{w}}_{i,0} := {\widetilde{w}_i(0)}  $) and introduce an approximation $u_{0,h} \in W_{h,0}$ to $u_0$.
  \item[ b)] Then, for given $n\geq 0$, compute the approximate state $u^n_h$ by solving
  \begin{equation} \label{3.10b}
\left\{
\begin{array}
[c]{ll}%
\vspace{.3cm}
u^{n,0}_h = {u}_{h,0} \, ,  \, \ \ u^{n,1}_{h} = u^{n,0}_h + (\Delta t) \cdot {u}_{h,1} \ , \\

\vspace{.2cm}
 \displaystyle \int_{\Omega}\biggl(\frac{1}{(\Delta t)^2}\bigl(u^{n, m+1}_h - 2 u^{n, m}_h +  u^{n, m-1}_h) z + \nabla u^{n, m+1}_h \cdot \nabla z  \biggr) \,dx  = 0 \, , \\  \vspace{.2cm}
  %\hspace{.7cm} = \displaystyle \int_{\Omega}\bigl( \hat{v}^n(x, t^{m+1})\mathds{1}_{\mathcal{O}_{1}}  \bigr)z \,dx,  \\
 \vspace{.2cm}
% \ \ \ \ \ \ \ \ \ \  \ \mbox{and} \ \ \ \ \ \ \\
\forall z \in W_{h,0}, \ \ u^{n, m+1}_h \in W_{h,0}, \ \ m=1,...,M-1,
\end{array}
\right.
   \end{equation}

and assuming that $u^{n, m}_h$ and $u_{2}(x, t^m)$ are known, compute the approximate adjoint states $p^{n, m}_{h}$ (for $u^{n, m}_h$),  by solving
  \begin{equation} \label{3.12c}
\left\{
\begin{array}
[c]{ll}%
\vspace{.3cm}
p^{n,M}_{h} = 0 \ ,  \ \ \ \ p^{n,M-1}_{h} = 0  \\
\vspace{.2cm}
\displaystyle \int_{\Omega}\biggl( \frac{1}{(\Delta t)^2}\bigl(p^{n,m+1}_{h} - 2p^{n,m}_{h} + p^{n,m-1}_{h})z + \nabla p^{n,m-1}_{h} \cdot \nabla z  \biggr) \,dx \\
= \displaystyle \int_{\Omega} (u^{n, m-1}_h - u_{2}(x, t^{m-1}))\cdot z \, dx \\
 \vspace{.2cm}
% \ \ \ \ \ \ \ \ \ \  \ \mbox{and} \ \ \ \ \ \ \\
\forall z \in W_{h,0}, \ \ p^{n,m}_{h} \in W_{h,0}, \  \ m= M-1, \, M-2, \, ... \, , \, 1
\end{array}
\right.
   \end{equation}

\item[c)]
Now, for given $n \geq 0$ consider known  $(\psi^{n}_h(0), \, \psi^{n}_{1,h}(0)) \in (W^{\Delta t}_{h,0}, \, W^{\Delta t}_{h,0})$ an approximation to $(\psi(0), \ \psi'(0)) \in W\times W$, and compute the approximate state $\psi^n_h$ to $\psi$, solving

 \begin{equation} \label{3.11b}
\left\{
\begin{array}
[c]{ll}%
\vspace{.3cm}
\psi^{n,0}_h = \psi^{1,0}_h = 0 \, ,  \, \ \ \psi^{n,1}_{h} = \psi^{n,0}_h + (\Delta t) \cdot {\psi}_{h,1} \ , \\

\vspace{.2cm}
 \displaystyle \int_{\Omega}\biggl(\frac{1}{(\Delta t)^2}\bigl(\psi^{n, m+1}_h - 2 \psi^{n, m}_h +  \psi^{n, m-1}_h) z + \nabla \psi^{n, m+1}_h \cdot \nabla z  \biggr) \,dx  = 0 \, , \\
  %\hspace{.7cm} = \displaystyle \int_{\Omega}\bigl( \hat{v}^n(x, t^{m+1})\mathds{1}_{\mathcal{O}_{1}}  \bigr)z \,dx,  \\
 \vspace{.2cm}
% \ \ \ \ \ \ \ \ \ \  \ \mbox{and} \ \ \ \ \ \ \\
\forall z \in W_{h,0}, \ \ \psi^{n, m+1}_h \in W_{h,0}, \ \ m=1,...,M-1,
\end{array}
\right.
   \end{equation}

in addition compute the approximate adjoint states $\varphi^{n, m}_{h},$ with $m = M-1, \, M-2, \, \ldots , \, 1$ (where $\varphi^{n, M}_{h}$ and $\varphi^{n, M-1}_{h}$ are known), by
\begin{equation} \label{3.12c}
\left\{
\begin{array}
[c]{ll}%
\vspace{.3cm}
\varphi^{n,M}_{h} = 0 \ ,  \ \ \ \ \varphi^{n,M-1}_{h} = 0  \\
\vspace{.2cm}
\displaystyle \int_{\Omega}\biggl( \frac{1}{(\Delta t)^2}\bigl(\varphi^{n,m+1}_{h} - 2\varphi^{n,m}_{h} + \varphi^{n,m-1}_{h})z + \nabla \varphi^{n,m-1}_{h} \cdot \nabla z  \biggr) \, dx \\
\vspace{.3cm}
= \displaystyle \int_{\Omega} \psi^{n, m-1}_{h} \cdot z \, dx \\

% \ \ \ \ \ \ \ \ \ \  \ \mbox{and} \ \ \ \ \ \ \\
\forall z \in W_{h,0}, \ \ \varphi^{n,m}_{h} \in W_{h,0}, \  \ m= M-1, \, M-2, \, ... \, , \, 1
\end{array}
\right.
   \end{equation}
and, finally set
\begin{equation}
%v^{n+1}_i = - \frac{1}{\mu} \Bigl[ \alpha \varphi^{\Delta t}_{1,h}  + (1 - \alpha) \varphi^{\Delta t}_{2,h} \Bigr] \bigg{|}_{\omega \times (0,T)} .
{\widetilde{w}_1}^{n+1} = - \varphi^{n}_x \bigg{|}_{\hat{\Sigma}_1}
\end{equation} \label{CONT2}
and
\begin{equation}
%v^{n+1}_i = - \frac{1}{\mu} \Bigl[ \alpha \varphi^{\Delta t}_{1,h}  + (1 - \alpha) \varphi^{\Delta t}_{2,h} \Bigr] \bigg{|}_{\omega \times (0,T)} .
{\widetilde{w}_2}^{n+1} = \frac{1}{{\sigma}} p^{n}_x \bigg{|}_{\hat{\Sigma}_2} ,
\end{equation} \label{CONT1}
with $\sigma$-fixed.
\end{description}

\section{Illustrative Numerical Examples}\label{Sec7}

\noindent

Thanks to the results obtained in the anterior sections and theoretical results obtained in \cite{Je1000}, we can consider for  each
$\displaystyle \widetilde{w}_1$, the \textit{Nash equilibrium} $\displaystyle \widetilde{w}_2$ associated to solution $\displaystyle u$
of (\ref{eq1.3.1}).  The computations have been performed using {\it Freefem{\tiny\bf ++}}, which is a high performance free software designed to solve problems of PDEs (see \cite{FH}).
As its name implies, it is a free software based on the Finite Element Method (more details are available at {{https://freefem.org/}}).  The graphic representations are obtained in combination with \textsc{MatLab}.
For all experiments the number of time steps in $M=100$ (that gives $\Delta t = T/M $). We consider $u_2 = 10$ fixed, the initial conditions  $u\Big{|}_{\Omega_0}$ are given by $u(0)=0$ and $u'(0)=0$. All initial and boundary conditions were programmed considering the information provided in system (\ref{eq3.139}).

We consider the interval $\widehat{\Sigma}_0 =(0,T)$ as control domain, where $\widehat{\Sigma}_1 = (T/2 , \, T)$ and $\widehat{\Sigma}_2 = (0 , \, T/2)$.
As the time for the problem must satisfy (\ref{hT}), with $k$ defined by (\ref{hT10}), we define $T_c$ as time of control (with $T_c > T$) and $k$ are fixed by $$ \displaystyle T_c = \frac{e^{\frac{2k(1+k)}{(1-k)^3}}}{k} \,   \ \ \ \mbox{and} \ \ \ \displaystyle k=\frac{1}{4}.$$

   Now, we present several tests for the algorithms in the section (\ref{Sec3.3}). Considering  $\varepsilon = 10^{-5}$ and the stopping criterion is determined by:
   %\begin{equation}\label{7.1a}
   $$
\frac{\|(\widetilde{w}^{n+1}_1 , \, \widetilde{w}^{n+1}_2) - (\widetilde{w}^n_1 , \widetilde{w}^n_2)\|}{\|(\widetilde{w}^{n+1}_1 , \widetilde{w}^{n+1}_2)\|} \leq \varepsilon.
   $$
   %\end{equation}

\begin{figure}[ht]
\centering
\includegraphics[width=0.7\linewidth]{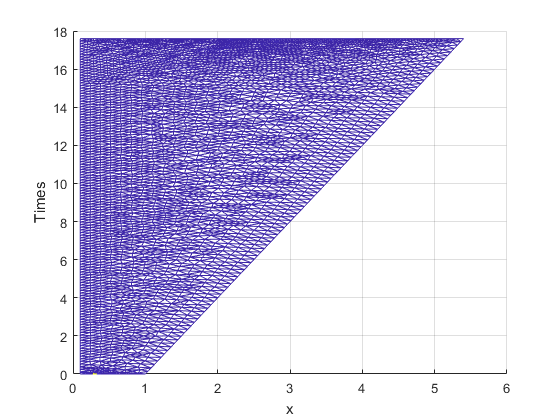} % leia abaixo
%\captionsetup{labelsep=none}
\vskip-.5cm
\caption{Domain for $k= 1/4$ and $T=2T_c$. Nb of vertices = 2916, Nb of triangles = 5526, Border length = 41.936 \ .}
\label{figura1}
\centering
\begin{multicols}{2}
\includegraphics[width=1.1\linewidth]{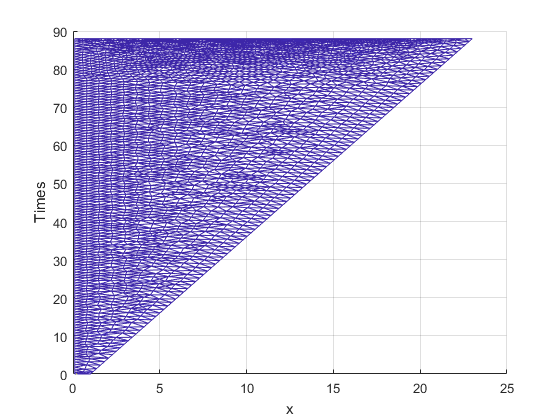}
\vskip-.4cm
         \caption{Domain for $k= 1/4$ and $T=5T_c$. Nb of vertices = 2319, Nb of triangles = 4332, Border length = 202.484 \ .}
\includegraphics[width=1.1\linewidth]{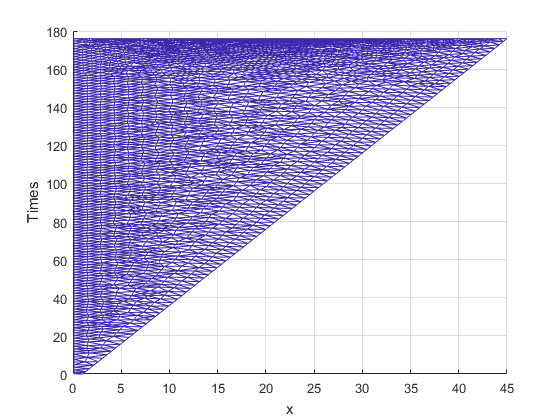}
\vskip-.4cm
   \caption{Domain for $k= 1/4$ and $T = 10 T_c$. Nb of vertices = 2236, Nb of triangles = 4166, , Border length = 403.167 \ .}
\end{multicols}
\label{figura2}
%{For all experiments, we consider the domain generated by $k= 1/4$.}
\end{figure}

\begin{figure}[htbp]
\centering
\includegraphics[width=0.7\linewidth]{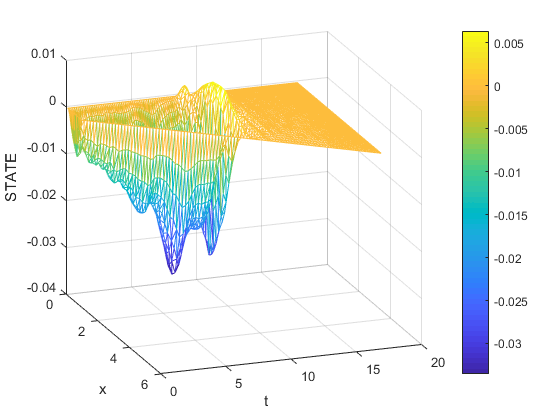} % leia abaixo
%\captionsetup{labelsep=none}
\vskip-.5cm
\caption{ Final state for the approximate solution $u^n$ in the time $T=T_c$. Iterates to the stopping criterion: 6.}
\label{figura1}
\begin{multicols}{2}
\centering
\includegraphics[width=1.1\linewidth]{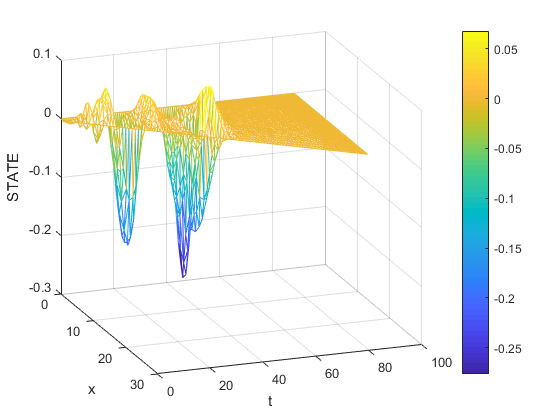}\\
\vskip-.4cm
\caption{ Final state for the approximate solution $u^n$ in the time $T=5 \cdot T_c$. Iterates to the stopping criterion  : 7.}
\includegraphics[width=1.1\linewidth]{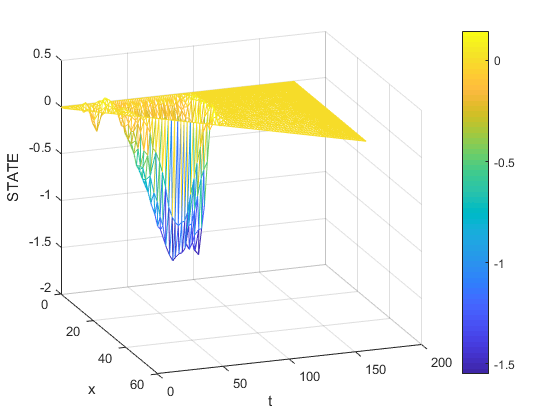}
\vskip-.4cm
\caption{  Final state for the approximate solution $u^n$ in the time $T= 10  T_c$. Iterates to the stopping criterion : 8.}
\end{multicols}
{Final states with change $T = T_c , \, 2 \cdot T_c, \, \ldots , 10 \cdot T_c \ $, fixed ${u}_2 = 10$ and $\sigma=10^2$. Maximum number of iterates = 100.}
\end{figure}
\label{figura3}

\begin{figure}[htbp]
\centering
\includegraphics[width=1.0\linewidth]{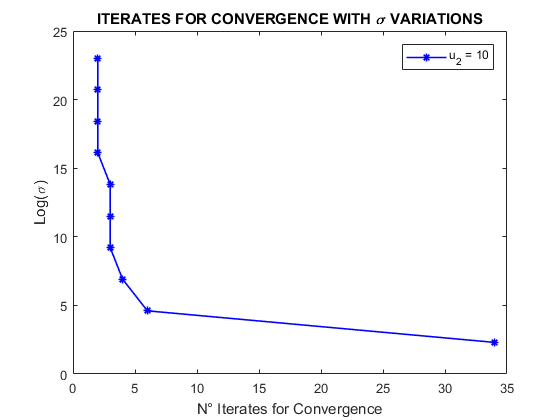} % leia abaixo
%\captionsetup{labelsep=none}
\vskip-.1cm
\caption{The number of iterations needed for convergence criterion when $\sigma=10, 10^2,  ... , 10^{10}$.}
\label{figura6}

\begin{multicols}{2}
\centering
\includegraphics[width=1.1\linewidth]{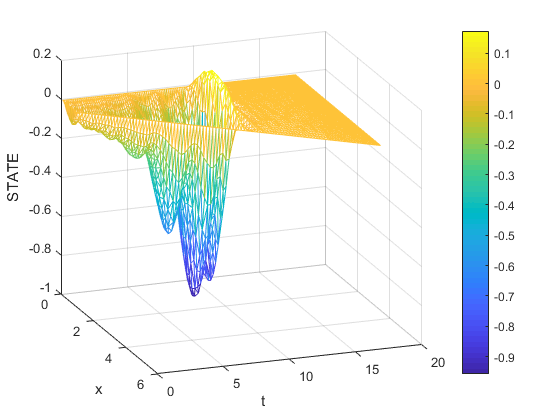}\\
         \centering\caption{Final state in $T=T_c$ and $\sigma=10$.}
\includegraphics[width=1.1\linewidth]{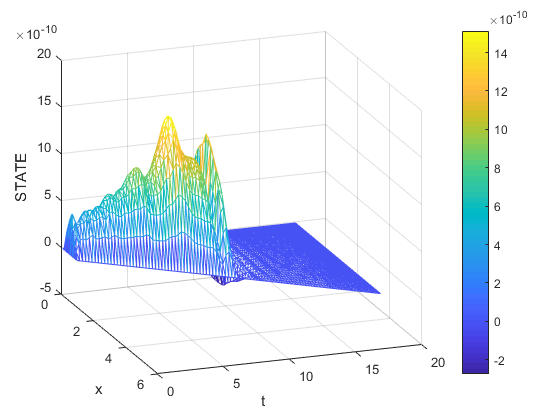}\\
         \centering\caption{Final state in $T=T_c$ and $\sigma=10^{10}$.}
\end{multicols}
%\vskip-.1cm
%\caption{A comparison of the influence of $\sigma$ in the final states plots, when $\sigma = 10$ and $\sigma = 10^{10}$, with $k=0.25$ and $T=Tc$.}
\end{figure}\label{mu-states}

\begin{table}[htbp]
%\centering
\begin{tabular}{llll}
%\begin{tabular}{c c c c}
\hline
\textsc{Times} & \hspace{1.2cm} \textsc{$ \displaystyle \| u - u^n \|_{L^2(\hat{Q})} < \varepsilon$} & \hspace{1.3cm} $\displaystyle \sum^2_{i=1} \| \tilde{w}_i - \tilde{w}_i^n \|_{L^2(\hat{\Sigma})} < \varepsilon$ & \hspace{.7cm} \textsc{Iterates} \tabularnewline
\hline \vspace{.05cm}
\hspace{.1cm}
$T_c$            &  \hspace{1.3cm} $4.23731 \cdot 10^{-6}$      & \hspace{1.6cm} $8.94841 \cdot 10^{-8}$    & \hspace{1.2cm} 6 \\ % \tabularnewline\hline
\vspace{.05cm}
$\ 2\cdot T_c$     &      \hspace{1.3cm} $1.39288 \cdot 10^{-5}$       & \hspace{1.6cm} $6.04906\cdot 10^{-6}$   & \hspace{1.2cm} 7 \\ %\tabularnewline\hline
\vspace{.05cm}
$\ 3 \cdot T_c$           & \hspace{1.3cm} $5.38065 \cdot 10^{-5}$       & \hspace{1.6cm} $3.47866 \cdot 10^{-7}$    & \hspace{1.2cm} 7 \\%\tabularnewline\hline
\vspace{.05cm}
$\ 4 \cdot  T_c$           & \hspace{1.3cm} $9.46165 \cdot 10^{-5}$       & \hspace{1.6cm} $7.61072 \cdot 10^{-7}$   & \hspace{1.2cm} 7 \\ % \tabularnewline\hline
\vspace{.05cm}
$\ 5 \cdot  T_c$           & \hspace{1.3cm} $9.84504 \cdot 10^{-5}$      & \hspace{1.6cm} $9.49046 \cdot 10^{-7}$    & \hspace{1.2cm} 7 \\ % \tabularnewline\hline
\vspace{.05cm}
$\ 6 \cdot T_c$           & \hspace{1.3cm} $2.22432 \cdot 10^{-5}$       &  \hspace{1.6cm} $1.62926 \cdot 10^{-7}$    & \hspace{1.2cm} 8 \\ %\tabularnewline\hline
\vspace{.05cm}
$\ 7 \cdot T_c$           & \hspace{1.3cm} $8.91691 \cdot 10^{-5}$       & \hspace{1.6cm} $6.78541 \cdot 10^{-7}$   & \hspace{1.2cm} 8 \\ %\tabularnewline\hline
\vspace{.05cm}
$\ 8 \cdot T_c$           & \hspace{1.3cm} $2.50240 \cdot 10^{-4}$      & \hspace{1.6cm} $1.41019 \cdot 10^{-6}$     & \hspace{1.2cm} 8 \\ % \tabularnewline\hline
\vspace{.05cm}
$\ 9 \cdot T_c$           & \hspace{1.3cm} $4.04092 \cdot 10^{-4}$       & \hspace{1.6cm} $1.85044 \cdot 10^{-6}$    & \hspace{1.2cm} 8 \\ % \tabularnewline\hline
\vspace{.05cm}
$10 \cdot T_c$           &\hspace{1.3cm} $5.41312 \cdot 10^{-4}$      &  \hspace{1.6cm} $2.14904 \cdot 10^{-6}$    & \hspace{1.2cm} 8 \tabularnewline\hline
%\hline
\end{tabular}\label{EX1}
\caption{ Domains construction in function of time $T=T_c$.  The maximum number of iterates = 100, $k=1/4$, $\sigma=10^2$  and ${u}_2 = 10$ fixed.}\label{tab0}
\end{table}

\begin{table}[htbp]
%\centering
\begin{tabular}{llll}
%\begin{tabular}{c c c c}
\hline
\textsc{Times} & \hspace{1.2cm} \textsc{$N.^{\circ}$ Vertices} & \hspace{1.4cm} \textsc{$N.^{\circ}$ Triangles} & \hspace{.8cm} \textsc{Border Length} \tabularnewline
\hline \vspace{.05cm}
\hspace{.1cm}
$T_c$            &  \hspace{1.6cm} $2916$      & \hspace{2.cm} $5526$    & \hspace{1.4cm} \ 41.936 \\ % \tabularnewline\hline
\vspace{.05cm}
$\ 2\cdot T_c$     &      \hspace{1.6cm} $2580$       & \hspace{2.cm} $4854$   & \hspace{1.4cm} \ 82.073 \\ %\tabularnewline\hline
\vspace{.05cm}
$\ 3 \cdot T_c$           & \hspace{1.6cm} $2411$       & \hspace{2.cm} $4516$    & \hspace{1.4cm} 122.210 \\%\tabularnewline\hline
\vspace{.05cm}
$\ 4 \cdot  T_c$           & \hspace{1.6cm} $2365$       & \hspace{2.cm} $4424$   & \hspace{1.4cm} 162.347 \\ % \tabularnewline\hline
\vspace{.05cm}
$\ 5 \cdot  T_c$           & \hspace{1.6cm} $2319$      & \hspace{2.cm} $4332$    & \hspace{1.4cm} 202.484 \\ % \tabularnewline\hline
\vspace{.05cm}
$\ 6 \cdot T_c$           & \hspace{1.6cm} $2309$       &  \hspace{2.cm} $4312$    & \hspace{1.4cm} 242.620 \\ %\tabularnewline\hline
\vspace{.05cm}
$\ 7 \cdot T_c$           & \hspace{1.6cm} $2316$       & \hspace{2.cm} $4326$   & \hspace{1.4cm} 282.757 \\ %\tabularnewline\hline
\vspace{.05cm}
$\ 8 \cdot T_c$           & \hspace{1.6cm} $2273$      & \hspace{2.cm} $4240$     & \hspace{1.4cm} 322.894 \\ % \tabularnewline\hline
\vspace{.05cm}
$\ 9 \cdot T_c$           & \hspace{1.6cm} $2246$       & \hspace{2.cm} $4186$    & \hspace{1.4cm} 363.031 \\ % \tabularnewline\hline
\vspace{.05cm}
$10 \cdot T_c$           &\hspace{1.6cm} $2236$      &  \hspace{2.cm} $4166$    & \hspace{1.4cm} 403.167 \tabularnewline\hline
%\hline
\end{tabular}\label{EX2}
\caption{ Domains construction in function of time $T=T_c$.  The maximum number of iterates = 100, $k=1/4$, $\sigma=10^2$  and ${u}_2 = 10$ fixed.}\label{tab0}
\end{table}



\begin{table}[htbp]
\begin{tabular}{lll}
\hline
\hspace{.1cm} ${\sigma}$ \hspace{2.8cm} & \textsc{Iterates} \hspace{2.8cm} & \textsc{{Error for stopping criteria}} \tabularnewline
\hline
 $10$             & \hspace{.3cm} $34$      &  \hspace{.8cm} $7.77814 \times 10^{-6}$  \\ %   \tabularnewline\hline
$10^2$           & \hspace{.4cm} $6$      &  \hspace{.8cm} $8.94080 \times 10^{-6}$  \\ %   \tabularnewline\hline
$10^3$           & \hspace{.4cm} $4$      & \hspace{.8cm} $2.00607 \times 10^{-6}$ \\ % \tabularnewline\hline
$10^4$           & \hspace{.4cm} $3$      & \hspace{.8cm}  $3.10847\times 10^{-6}$  \\ %  \tabularnewline\hline
$10^5$           & \hspace{.4cm} $3$      & \hspace{.8cm}  $3.10853 \times 10^{-8}$  \\ % \tabularnewline\hline
$10^6$           & \hspace{.4cm} $3$      & \hspace{.8cm} $3.10861 \times 10^{-10}$  \\% \tabularnewline\hline
$10^7$           & \hspace{.4cm} $2$      & \hspace{.8cm} $5.61446  \times 10^{-6}$ \\
$10^8$           & \hspace{.4cm} $2$      & \hspace{.8cm} $5.61446 \times 10^{-7}$ \\
$10^9$           & \hspace{.4cm} $2$      & \hspace{.8cm} $5.61446  \times 10^{-8}$ \\
$10^{10}$           & \hspace{.4cm} $2$      & \hspace{.8cm} $2.61444\times 10^{-8}$
  \tabularnewline\hline
%\hline
\end{tabular}
\caption{ The iterates and convergence error with change ${\sigma}$. The maximum number of iterates = 100, $k=1/4$, $T=T_c$  and ${u}_2=10$ fixed.}\label{tab1}

\begin{remark}{
When $k$ increases in (0,1), with $T=T_c$ and $\sigma = 10^2$ fixed, the convergence of the algorithm does not show good results. But, considering increasing $\sigma$ $(\sigma = 10^3, 10^4, 10^5, ...)$, good convergence results are obtained.}
\end{remark}
\end{table}

\section{Some additional comments and conclusions}\label{Sec8}
%\noindent
We have presented a numerical approach for the hierarchical control problem to the wave equation, with the mobile boundary and the controls acting on an  piece of the border. We use the results proven in \cite{Je1000} to ensure the validity of the results that underlie the numerical part developed.  In the numerical part, we use a combination of tools: Finite Element Method (in space) and Finite Difference (in time), adding a  fixed point algorithm to evaluate the computational convergence of the obtained results. We have established also the feasibility of simulating the problems  on which hierarchical control acts in the moving boundary.

These results can help numerical and computational advances in other types of equilibrium problems with controls acting on the moving limit, such as \textit{Stackelberg-Pareto}, \textit{Pareto}, \textit{Nash}, among others. It can also be extended into similar analyses for other types of hierarchical control problems, such as \textit{Heat equation, Stokes, Navier-Stokes, Schr\"odinger} (in \cite{Da} some results are presented), among others.

%\begin{acknowledgements}
%\end{acknowledgements}

\section*{Acknowledgements}

This work was supported by Universidade Estadual do Piau\'i-UESPI and EDITAL FAPEPI/MCT/CNPq N$^{o}$ 007/2018:Programa de Infraestrutura para Jovens Pesquisadores/ Programa Primeiro Projetos (PPP).

%The authors thank the support of Universidade Estadual do Piau\'i. This work was supported by the ''EDITAL FAPEPI/MCT/CNPq N$^{o}$ 007/2018: Programa de Infraestrutura para Jovens Pesquisadores/ Programa Primeiro Proje- tos (PPP)``.

\section*{Bibliography}

\end{document}